\documentclass[10pt, english,oneside,reqno]{amsart}
\usepackage[utf8]{inputenc}
\usepackage[english]{babel}
\usepackage{hyperref}
\usepackage{color}
\usepackage{enumerate}
\usepackage{fancyhdr}

\usepackage{amsmath}
\usepackage{amscd}
\usepackage{amsfonts}
\usepackage{amsthm}
\usepackage{amstext}
\usepackage{amssymb}
\usepackage{amsbsy}
\usepackage{mathrsfs}
\usepackage{tikz-cd}
\usepackage{mathtools}
\usepackage{enumitem}
\newtheorem*{corollary*}{Corollary}

\usepackage{float}
\usepackage{caption}
\usepackage{subcaption}
\usepackage{graphicx}

\usepackage[alphabetic,abbrev]{amsrefs}

\usepackage{url}

\theoremstyle{plain}
\newtheorem{theorem}{Theorem}[section]

\newtheorem{proposition}[theorem]{Proposition}
\newtheorem{corollary}[theorem]{Corollary}
\newtheorem{lemma}[theorem]{Lemma}

\theoremstyle{remark}

\numberwithin{equation}{section}
\numberwithin{figure}{section}
\definecolor{esperance}{rgb}{0.1,0.5,0.3}

\newcommand{\E}{\mathbb{E}}
\newcommand{\Prob}{\mathbb{P}}
\newcommand{\R}{\mathbb R}

\DeclareMathOperator{\supp}{supp}
\DeclareMathOperator{\Var}{Var}

\DeclareMathOperator{\diam}{diam}
\hypersetup{
  colorlinks=false,
  pdfborder={0 0 1},
  linkbordercolor={0 1 0},
  citebordercolor={0 1 0},
  urlbordercolor={0 1 0}
}

\makeatletter
\renewcommand{\eqref}[1]{%
  \textup{\tagform@{%
  \begingroup
  \hypersetup{linkbordercolor={1 0 0}}%
  \hyperref[#1]{\ref*{#1}}%
  \endgroup
  }}%
}
\makeatother

\DeclareRobustCommand{\thmref}[1]{%
  \begingroup
  \hypersetup{linkbordercolor=red}%
  \hyperref[#1]{\ref*{#1}}%
  \endgroup
}
\title[On $L^p$ density of self similar measures and Bernoulli Convolutions]{On $L^p$ density of self similar measures and Bernoulli Convolutions}
\date{}

\begin{document}

\begin{abstract}
We prove that an absolutely continuous self-similar measure has an $L^p$ density for some $p>1$. 
Moreover, explicit examples of Bernoulli convolutions with a continuous density are established.
\end{abstract}

\author[sk]{Samuel Kittle}
\author[ck]{Constantin Kogler}

\email{s.kittle@ucl.ac.uk, kogler@ias.edu}

    \address{Samuel Kittle, Department of Mathematics, University College London, 25 Gordon Street, London WC1H 0AY, United Kingdom}

    \address{Constantin Kogler, Institute for Advanced Study, 1 Einstein Dr, Princeton, NJ 08540, United States of America}
\maketitle

\section{Introduction}

Self-similar measures are central objects in fractal geometry \cites{Fal14,
BSS23}. Understanding when they are absolutely continuous and establishing
results about their regularity is intricate. In this paper, we establish the
fundamental property that an absolutely continuous self-similar measure has a
density in $L^p$ for some $p > 1$. Moreover, we apply the quantitative proof method of the
latter result to establish explicit examples of Bernoulli convolutions with a
continuous density going beyond the only previously known case of Garsia
numbers. All new results of this paper are  formalized in Lean.

A contracting similarity on $\mathbb{R}^d$ is a map of the form $x \mapsto \rho
U x + b$ with $\rho \in (0,1)$, $U$ a $d\times d$ orthogonal matrix and $b \in
\mathbb{R}^d$. Let $\mu=\sum_{i=1}^k p_i\delta_{g_i}$ be a probability measure
supported on finitely many contracting similarities on $\R^d$. The self-similar
measure of $\mu$ is the unique probability measure $\nu$ on $\mathbb{R}^d$
satisfying
\begin{equation}\label{eq:stationarity}
 \mu*\nu:=\sum_{i=1}^k p_i(g_i)_*\nu=\nu.
\end{equation}
Ever since Hutchinson's 1981 paper \cite{Hut81}, where it is proved that they
exist and have compact support, self-similar measures have been well studied.
We establish the following observation about self-similar measures.

\begin{theorem}\label{thm:main}
Let $\mu$ be a probability measure supported on finitely many contracting
similarities of $\R^d$ for $d\geq 1$ and denote by $\nu$ the self-similar measure of $\mu$. If
$\nu$ is absolutely continuous with respect to the Lebesgue measure, then the
density of $\nu$ is in $L^p(\R^d)$ for some $p>1$.
\end{theorem}

The proof of Theorem~\thmref{thm:main} is short and combines self-similarity 
with a general higher integrability result. The key lemma of the proof
of Theorem~\thmref{thm:main} and of this paper in general is 
Lemma~\ref{lem:weak}, which is a version of the reverse H\"older criterion 
of Spadaro \cite{Spa12}, as further developed by Kinnunen--Kurki--Mudarra \cite{KKM22}. 
Lemma~\ref{lem:weak} assumes for an absolutely continuous measure $\nu$ on $\mathbb{R}^d$ 
control over the mass of $\nu(E)$ for sets $E$ of small measure and then establishes a 
regularity gain for $\nu$. The proof of Theorem~\ref{thm:main} will then follow by 
showing that self-similarity prevents the mass of $\nu(E)$ being too large.
Similar higher integrability arguments are used by Giaquinta--Modica
\cite{GM86}*{Theorems~3.1 and~4.1} in the regularity theory of elliptic
equations and by De Lellis--Spadaro \cite{DLS14}*{Section~5} in
geometric measure theory.

The exponent $p>1$ in Theorem~\thmref{thm:main} and Lemma~\ref{lem:weak}
can be estimated quantitatively from the mass that $\nu$ assigns to sets of
small Lebesgue measure. This insight allows us to use recent methods for establishing 
explicit examples of absolutely continuous
self-similar measures \cites{Var19, Kit24, KK25} to obtain suitable measure
estimates, yielding quantitative $L^p$ results that we use to construct
Bernoulli convolutions with continuous densities. 

The Bernoulli convolution $\nu_{\lambda}$ of parameter $\lambda \in (0,1)$ is
the self-similar measure of the balanced probability measure on the similarities
$x\mapsto \lambda x + 1$ and $x \mapsto \lambda x - 1$ with $x \in \mathbb{R}$.
Another way to describe them is as the law of the random variable
$\sum_{n=0}^{\infty}\xi_n\lambda^n$, where the $\xi_n$ are independent and satisfy
$\Prob(\xi_n=1)=\Prob(\xi_n=-1)=1/2$ for all $n\geq 0$. The study of Bernoulli convolutions has a
rich history and we refer to surveys by Peres--Schlag--Solomyak \cite{PSS00} and
Varj\'u \cite{Var18}.

We highlight that Solomyak \cite{Sol95} showed in 1995 that $\nu_{\lambda}$
is in $L^2$ for almost all $\lambda \in (1/2,1)$. Another important result is
due to Shmerkin \cite{Shm19}*{Theorem~1.3},
who showed that $\nu_\lambda$ has a density in $L^p$ for every finite 
$p>1$ and all parameters $\lambda\in(1/2,1)$ outside a set of Hausdorff dimension zero. 

Splitting the series $\sum_{n\ge0}\xi_n\lambda^n$ into even and odd $n$ gives
$\nu_\lambda=\nu_{\lambda^2}*(D_\lambda)_*\nu_{\lambda^2}$ with
$D_\lambda(x)=\lambda x.$ Therefore, by \cite{Shm19}*{Theorem~1.3} applied to $L^2$, 
for all $\lambda\in(1/\sqrt2,1)$ except for a set of Hausdorff dimension zero, the factors $\nu_{\lambda^2}$ and
$(D_\lambda)_*\nu_{\lambda^2}$ are both in $L^2$. As the 
convolution of two $L^2$ functions is continuous, it follows that
$\nu_\lambda$ has a continuous density for $\lambda \in
(1/\sqrt2,1)$ outside a set of Hausdorff dimension zero. By a similar argument using Young's inequality, it follows
that if $\nu_{\lambda} \in L^p$, then $\nu_{\lambda^{1/m}}$ has a continuous
density for an integer $m$ satisfying $m \geq p/(p-1)$. It is therefore a direct
consequence of Theorem~\thmref{thm:main} that if $\nu_{\lambda}$ is absolutely
continuous, then $\nu_{\lambda^{1/m}}$ has a continuous density for
sufficiently large $m$. As we are moreover able to make $p$ quantitative, we are
led to fully explicit examples of Bernoulli convolutions with a continuous density.

Recall that the Mahler measure $M_{\lambda}$ of an algebraic number $\lambda$ 
is defined as $$M_{\lambda} = |a|\prod_{i} \max\{1, |\lambda_i| \},$$ where the 
primitive minimal polynomial of $\lambda$ in $\mathbb{Z}[x]$ has the form 
$a(x - \lambda_1)\cdots (x- \lambda_d)$.
Garsia \cite{Gar62} proved in 1962 that $\nu_{\lambda}$ has an $L^{\infty}$ density if
$\lambda \in (1/2,1)$ is a Garsia number, that is, $\lambda$ satisfies $M_{\lambda} = 2$.
There was no progress on understanding absolute continuity for explicit
$\lambda$ until Varjú \cite{Var19} proved in 2019 that $\nu_{\lambda}$ is
absolutely continuous if $\lambda > 1 - c\min\{\log M_{\lambda}, (\log
M_{\lambda})^{-2}\}$ for $c > 0$ an absolute constant. Varjú's result was
improved in \cite{KK25}, addressing general self-similar measures, by showing
that $\nu_{\lambda}$ is absolutely continuous if $\lambda > 1 - c\min\{\log
M_{\lambda}, (\log \log M_{\lambda})^{-2}\}$ for $c > 0$ an absolute constant.
Yu \cite{Yu22}  proved in 2022 that $\nu_{\lambda}$ has a continuous density if
$\lambda$ is a Garsia number satisfying $\lambda\in [1/\sqrt{2},1)$.

In this paper, we give the first explicit non-Garsia examples of Bernoulli
convolutions with a continuous density.

\begin{theorem}\label{thm:algebraic}
If $\lambda\in(1/2,1)$ is algebraic of degree $d_{\lambda}$, has Mahler measure
$M_{\lambda}$ and satisfies
\begin{equation}\label{eq:main-condition}
 \lambda\ge1-\left(\frac{\min\{1,\log M_\lambda\}}
 {1+d_\lambda+\log M_\lambda}\right)^{100},
\end{equation}
then $\nu_\lambda$ is absolutely continuous with a density in $C_c(\R)$.
\end{theorem}

We further comment on the proof of Theorem~\ref{thm:algebraic} in Section 3. 
To deduce Theorem~\ref{thm:algebraic} from Lemma~\ref{lem:weak}, a 
simplified version of the variance summation method from \cite{Kit25} and \cite{KK25}
is applied. The reason that, in contrast to \cites{Var19, KK25},
the condition from  Theorem~\ref{thm:algebraic} depends on the degree of $\lambda$
is because to apply Lemma~\ref{lem:weak} one needs
to control the measure from a fixed scale onwards. 
This requires explicit separation bounds, which we establish in Lemma~\ref{lem:algebraic-input} 
and compare to previous separation estimates immediately after the statement.

The exponent $100$ in Theorem~\ref{thm:algebraic} is chosen to make the argument elegant. In fact, a refinement 
of the argument presented in this paper to lower the exponent from $100$ to $3$ is available at \cite{KKsupp}.

\begin{theorem}\label{thm:algebraic_strong}(\cite{KKsupp})
If $\lambda\in(1/2,1)$ is algebraic of degree $d_{\lambda}$, has Mahler measure
$M_{\lambda}$ and satisfies
\begin{equation}\label{eq:main-condition_strong}
 \lambda\ge1-\left(\frac{\min\{1,\log M_\lambda\}}
 {1+d_\lambda+\log M_\lambda}\right)^{3},
\end{equation}
then $\nu_{\lambda^2}$ is absolutely continuous with a density in
$L^2(\mathbb{R})$ and therefore $\nu_{\lambda}$ has a density in $C_c(\R)$. 
\end{theorem}

When one applies Theorem~\thmref{thm:algebraic_strong} to rational numbers, one obtains
the following.

\begin{corollary}\label{cor:rational}
If $\lambda=1-p/q\in(1/2,1)$ is in lowest terms and
\[
 p\leq \frac{q}{(2 + \log q)^{3}},
\]
then $\nu_\lambda$ is absolutely continuous with a density in $C_c(\R)$. In particular, $\nu_{1 - \frac{1}{n}}$ 
is absolutely continuous with a density in $C_c(\R)$ for all $n\geq 588$.
\end{corollary}

\begin{proof}
As $d_\lambda=1$ and $M_\lambda=q\ge3$ the claimed inequality follows from
Theorem~\thmref{thm:algebraic_strong}. 
The desired result for $\nu_{1 - \frac{1}{n}}$ follows from a brief calculation.
\end{proof}

\subsection*{Machine Contribution Statement}
The main idea of this paper is due to GPT-6 Astra. On September 10th
2026, the second-named author (C.K.) asked GPT-6 Astra to find explicit examples
of self-similar measures with an $L^p$ density resulting in the discovery of
Theorem~\thmref{thm:main}. Various simplifications led to Lemma~\ref{lem:weak}.
Realising that Theorem 1.1 is quantitative, numerous conversations
with GPT-6 Astra led to
Theorem 1.2. The paper is written by the authors and we take responsibility for its correctness.

\subsection*{Lean Formalization}  GPT-6 Astra formalized all results from this paper 
and all necessary results from previous work, including Theorem~\ref{thm:algebraic_strong}. 
The Lean source code is available in the accompanying repository \cite{KogLean}.

\subsection*{Acknowledgments} The first-named author (S.K.) gratefully acknowledges support 
from the Heilbronn Institute for Mathematical Research. 
The second-named author (C.K.) holds a Postdoc Mobility Fellowship from the 
Swiss National Science Foundation (grant number
235409) and thanks the Institute for Advanced Study. 
The second-named author (C.K.) thanks Elon Lindenstrauss for his comments and 
interest in this work.

\section{Proof of Theorem~\thmref{thm:main}}

Write $|E|$ for the Lebesgue measure of a subset $E\subset \mathbb{R}^d$. All balls considered are open. If $B = B(x,r)$ is an open ball around $x \in \mathbb{R}^d$ and $r > 0$, then we denote for $a > 0$ by
$$aB = B(x,ar).$$ 

The main idea of the paper is to apply the following lemma to self-similar measures. It is a form of the weak $A_\infty$ criterion of Spadaro
\cite{Spa12}*{Theorem~1.1}, \cite{KKM22}*{Theorem~3.1 and Section~4} or \cite{AHT17}*{Corollary 5.5}.

\begin{lemma}[Weak reverse H\"older criterion]\label{lem:weak}
For every $d \geq 1$ write $C = 18^d$. Let $f\ge0$ belong to $L^1(\R^d)$. Suppose that for some $\varepsilon \in (0,1)$ and $\eta \in (0,C^{-1})$ it holds that
\begin{equation}\label{eq:weak-ainfty}
 |E|\le\varepsilon|B|
 \quad\Longrightarrow\quad
 \int_E f\,dx\le\eta\int_{2B}f\,dx
\end{equation}
for every ball $B$ of radius less than one and every measurable
$E\subset B$. Then $f\in L^p(\R^d)$ for $p = 1 + \frac{\log (C\eta)^{-1}}{\log(C\varepsilon^{-1})} > 1$.
\end{lemma}

\begin{proof}[Proof of Lemma~\ref{lem:weak}]
Write
\[
    A=C\varepsilon^{-1},
    \qquad
    \gamma=2\cdot6^d\eta \quad \text{ and } \quad
    F(t)=\int_{\{f>t\}}f\,dx.
\]
We claim that, for all sufficiently large $t$,
\begin{equation}\label{eq:lemma21-tail}
    F(At)\le\gamma F(t).
\end{equation}

Let $E=\{f>At\}$. By the Lebesgue differentiation theorem, almost all $x \in E$ satisfy that $\lim_{r \to 0} \frac{1}{|B(x,r)|} \int_{B(x,r)} f \, dx = f(x) > At > 2t$. For each such point $x\in E$, choose the ball $B_x$ centered at $x$ with largest radius whose average of $f$ is at least $2t$.
Continuity of ball averages gives
\[
    \frac{1}{|B_x|}\int_{B_x}f\,dx=2t
    \quad \text{ and } \quad
    \frac{1}{|aB_x|}\int_{aB_x}f\,dx\le2t
    \quad \text{ where } \quad a>1.
\]
For sufficiently large $t$, every ball $3B_x$ has radius less than one,
since $2t|B_x|\le\|f\|_1$.

Cover a compact subset $K\subset E$ of these Lebesgue points by finitely many
such balls and then choose a selection of disjoint balls $B_i$ with
$K\subset\bigcup_i3B_i$. Since for each $i$,
\[
    |E\cap3B_i|
    \le\frac{1}{At}\int_{3B_i}f\,dx
    \le\frac{2}{A}|3B_i|
    \le\varepsilon|3B_i|,
\]
the hypothesis \eqref{eq:weak-ainfty} yields
\[
    \int_K f\,dx
    \le\eta\sum_i\int_{6B_i}f\,dx
    \le2\cdot6^d\eta\,t\sum_i|B_i| = \gamma t\sum_i|B_i|.
\]
On the other hand, $2t|B_i|=\int_{B_i}f\,dx
    \le t|B_i|+\int_{B_i\cap\{f>t\}}f\,dx$.
Disjointness therefore gives $t\sum_i|B_i|\le F(t)$.
By inner regularity, $F(At)\le\gamma F(t)$, proving
\eqref{eq:lemma21-tail}.

Choose $T>0$ so that \eqref{eq:lemma21-tail} holds for $t\ge T$, and write $I_R=\int_T^R t^{p-2}F(t)\,dt$. Note that $\gamma A^{p-1}
    =\frac{2\cdot6^d\eta}{C\eta}
    =2\left(\frac13\right)^d\leq \frac23$. Therefore, for $R\ge AT$, changing variables in \eqref{eq:lemma21-tail} gives $I_R
    \le I_{AT}+\gamma A^{p-1}I_{R/A}
    \le I_{AT}+\frac23 I_R$.
Thus $I_R \leq 3I_{AT}$ and hence $I_R$ is uniformly bounded. Letting $R\to\infty$ and using
$F\le\|f\|_1$, Tonelli's theorem gives
\[
    \frac{1}{p-1}\int_{\mathbb R^d}f^p\,dx =\int_{\mathbb{R}^d} f(x) \left( \int_0^{f(x)} t^{p-2}\,dt \right)\, dx
    =\int_0^\infty t^{p-2}F(t)\,dt
    <\infty.
\] This proves that $f \in L^p(\mathbb{R}^d)$ and so concludes the proof.
\end{proof}

Theorem~\thmref{thm:main} follows by applying Lemma~\ref{lem:weak} together with self-similarity.

\begin{proof}[Proof of Theorem~\ref{thm:main}]
Without loss of generality we may normalize $\nu$ such that $\diam K=1$, where $K=\supp\nu$. For a similarity $g$, we write $\rho(g)$ for the contraction factor of $g$ and denote
\[
 \rho_* =\min_i\rho(g_i)>0
 \quad \text{ and } \quad
 \omega(t)=\sup\{\nu(E):|E|\le t\}.
\]
Absolute continuity gives $\omega(t)\to0$ as $t\downarrow0$.
We claim that for every ball $B=B(x,r)$ with $0<r<1$ and every
measurable $E\subset B$,
\begin{equation}\label{eq:local}
 \nu(E)\le\omega\!\left(c\frac{|E|}{|B|}\right)\nu(2B)
 \qquad  \text{for} \qquad c=|B(0,1)|\rho_*^{-d}.
\end{equation}
Since $\omega(t)\to0$, by \eqref{eq:local} it follows that the 
density function of $\nu$ verifies the hypothesis of 
Lemma~\ref{lem:weak} for every $\eta > 0$, which proves the theorem.

So it remains to prove \eqref{eq:local}. Let $\gamma_1, \gamma_2, \cdots$ be independent samples from $\mu$. Write $q_n = \gamma_1\cdots\gamma_n$ and consider the stopping time $\tau_{r} = \min\{ n \geq 1 \,:\, \rho(q_n) \leq r  \}$ and note $\rho_{*}r <\rho(q_{\tau_r}) \leq r$. It is well known that $q_{\tau_r} \nu$ has the same distribution as $\nu$, as is shown in \cite{Kit25}*{Lemma~2.19}, in other words $\E\bigl[(q_{\tau_r})_*\nu\bigr]=\nu.$ 

Consider a set $E \subset B$. Observe that $|q_{\tau_r}^{-1}E|=\rho(q_{\tau_r})^{-d}|E| \le c\frac{|E|}{|B|}$. Moreover, for $\nu(q_{\tau_r}^{-1}E) \neq 0$ it must hold that $q_{\tau_r}K\cap B \neq \emptyset$. So consider the event $\mathcal H=\{q_{\tau_r}K\cap B\ne\varnothing\}$.
On $\mathcal H$, the set $q_{\tau_r}K$ has diameter at most $r$, so it is contained
in $2B$ and therefore $\Prob(\mathcal H)\le\E\bigl[\nu(q_{\tau_r}^{-1}(2B))\bigr]
 =\nu(2B).$
Thus, it follows by the definition of $\omega$ that
\[
 \nu(E)=\E\bigl[\nu(q_{\tau_r}^{-1}E)\bigr]
 \le\omega\!\left(c\frac{|E|}{|B|}\right)\Prob(\mathcal H)
 \le\omega\!\left(c\frac{|E|}{|B|}\right)\nu(2B),
\]
which proves \eqref{eq:local}.
\end{proof}

\section{Bernoulli convolutions: Proof of Theorem~\ref{thm:algebraic}}

Throughout this section, write $d_{\lambda}$ for the degree 
of an algebraic number $\lambda\in(0,1)$ and
\[
 \chi_\lambda=-\log\lambda \qquad \text{ and } \qquad
 B_\lambda=\frac{1+d_\lambda+\log M_\lambda}
 {\min\{1,\log M_\lambda\}}.
\]
Thus \eqref{eq:main-condition} is equivalent to
$(1-\lambda)B_\lambda^{100}\le1$.
All logarithms are natural. The main result of this section is 
the following proposition.

\begin{proposition}[A quantitative gain over $L^1$]\label{prop:small-gain}
If $\lambda\in(1/2,1)$ is algebraic and
$\chi_\lambda B_\lambda^{10}\le10^{-100}$, then $\nu_\lambda$ is absolutely continuous with density $f_{\lambda}$ and
\begin{equation}\label{eq:small-gain}
 f_\lambda\in L^{1+\delta_\lambda}(\R)\quad \text{ for } \quad
 \delta_\lambda=\frac{1}{10(1-\log(1-\lambda))}.
\end{equation}
\end{proposition}

In Section~\ref{sec:proofmain} we deduce
Theorem~\ref{thm:algebraic} from Proposition~\ref{prop:small-gain} by applying Young's 
inequality and then prove Proposition~\ref{prop:small-gain} in Section~\ref{sec:detail} and
Section~\ref{sec:effective-detail-proof}.
The proof of Proposition~\ref{prop:small-gain} uses the detail of a measure, 
which is a concept
introduced by the first-named author \cite{Kit24}, the notion of 
entropy at a scale from Varjú \cite{Var19}
and a simplified version of
the variance summation method from \cite{Kit25} and \cite{KK25}.
To apply Lemma~\ref{lem:weak} to deduce Proposition~\ref{prop:small-gain}, we need to 
control $\nu_{\lambda}$ at a fixed explicit scale for which we use the 
explicit separation bound from Lemma~\ref{lem:algebraic-input}.

As the proof of Proposition~\ref{prop:small-gain} is an illuminating
series of ideas, we will not further outline it and 
simply let it unfold as presented below. We refer the reader to \cite{KK25} for
an exposition of the variance summation method.

\subsection{Deduction of Theorem~\ref{thm:algebraic}}\label{sec:proofmain}

We will use Young's inequality in the following form.

\begin{lemma}[Young's inequality for scaled convolutions]\label{lem:young-convolution}
Let $k\ge2$ be an integer, let $a_1,\ldots,a_k>0$, and let $\nu$ be a compactly supported
absolutely continuous probability measure on $\mathbb{R}$ with a density in $L^p(\R)$ for some
$p\ge k/(k-1)$. Then
\[
 \mathop{*}_{j=1}^{k}(D_{a_j})_*\nu
\]
has a density in $C_c(\mathbb{R})$, where $D_a(x)=ax$.
\end{lemma}

\begin{proof}
For $1\le j\le k$, let $f_j$ be the density of $(D_{a_j})_*\nu$. Compact support and
$p\ge k/(k-1)$ imply $f_j\in L^{k/(k-1)}$. By iterating Young's convolution inequality,
the function $g = f_1*\cdots*f_{k-1}$ is in $L^k$ since $\frac{k-1}{k/(k-1)}-(k-2)=\frac1k.$
Write $h=g*f_k$. H\"older's inequality
gives $\|h\|_\infty\le\|g\|_k\|f_k\|_{k/(k-1)}$ and, for every
$x,t\in\mathbb R$,
\[
 |h(x+t)-h(x)|
 \le \|g\|_k\,
 \|f_k(\,\cdot+t)-f_k\|_{k/(k-1)}.
\]
The right-hand side tends to zero as $t\to0$, by continuity of
translations in $L^{k/(k-1)}$. Hence $h$ is uniformly continuous.
Moreover, its support is contained in the sum of the compact
supports of the factors, so $h$ is compactly supported.
\end{proof}

The following basic estimate on $B_{\lambda}$ is used as well.

\begin{lemma}[Arithmetic input]\label{lem:easy_algebraic-input}
Let $\lambda\in(1/2,1)$ be an algebraic number. Then
\begin{equation}\label{eq:arithmetic}
 B_\lambda\ge3\quad \text{ and } \quad
 \log\frac1{1-\lambda}\le B_\lambda.
\end{equation}
\end{lemma}

\begin{proof}
    The first assertion in \eqref{eq:arithmetic} is immediate from the definition of
$B_\lambda$, since $B_\lambda\ge2+d_\lambda\ge3$.
Let $P_\lambda(T)=a\prod_{i=1}^{d_\lambda}(T-\lambda_i)$ be the primitive
integer minimal polynomial of $\lambda$, with $a>0$ and $\lambda_1=\lambda$.
Since $P_\lambda(1)$ is a nonzero integer and as 
$|1 - \lambda_i| \leq 1 + |\lambda_i| \leq 2 \max\{1, |\lambda_i| \}$ it follows that
$1\leq |P_\lambda(1)|\le2^{d_\lambda-1}M_\lambda(1-\lambda)$. Therefore, we have
\[
 \log\frac1{1-\lambda}
 \le(d_\lambda-1)\log2+\log M_\lambda\le B_\lambda.
\]
This proves \eqref{eq:arithmetic}.
\end{proof}

We now deduce Theorem~\thmref{thm:algebraic} from Proposition~\ref{prop:small-gain},
Lemma~\ref{lem:young-convolution} and Lemma~\ref{lem:easy_algebraic-input}.

\begin{proof}[Proof of Theorem~\thmref{thm:algebraic}]
We let $\theta = \lambda^k$ for some positive integer $k$ we will choose later. 
We wish to show that $\theta$ satisfies the conditions of Proposition \ref{prop:small-gain} 
and for Lemma~\ref{lem:young-convolution} that 
\begin{equation}
\delta_{\theta} \geq \frac{1}{k-1}. \label{eq:suff_continous}
\end{equation}

First note that since $\theta \leq \lambda$ we have $\delta_{\theta} \geq \delta_{\lambda}$. 
Therefore for \eqref{eq:suff_continous} to hold it is sufficient that 
$k \geq \frac{1}{\delta_{\lambda}} + 1 = 11 - 10 \log(1 - \lambda)$. 
Thus, we take $$k = \lfloor 20 - 20 \log(1 - \lambda) \rfloor$$ and claim that
\begin{equation}\label{eq:theta-small-gain}
 \chi_{\theta} B_{\theta}^{10} \leq 10^{-100}.
\end{equation}
As $B_\theta\ge1$, \eqref{eq:theta-small-gain} implies $\theta\in(1/2,1)$ and so
Proposition~\ref{prop:small-gain} gives $f_\theta\in L^{1+\delta_\theta}$.
Splitting the defining series modulo $k$ gives
$\nu_\lambda=\mathop{*}_{j=0}^{k-1}(D_{\lambda^j})_*\nu_\theta$, so
\eqref{eq:suff_continous} and Lemma~\ref{lem:young-convolution} give a continuous
density for $\nu_\lambda$, proving the theorem.

To prove \eqref{eq:theta-small-gain}, note that standard height identities
\cite{BG06}*{Section~1.5 and Section~1.6} give
\[
 \frac{\log M_\theta}{d_\theta}
 =k\frac{\log M_\lambda}{d_\lambda}
 \quad \text{ and } \quad \frac{d_\lambda}{k}\le d_\theta\le d_\lambda.
\]
Hence $\log M_\lambda\le\log M_\theta\le k\log M_\lambda$ and
$B_\theta\le kB_\lambda$. Write $t=-\log(1-\lambda)$.
Assumption~\eqref{eq:main-condition} gives $e^{-t}=1-\lambda\le B_\lambda^{-100}$ and therefore
\[
 t\ge100\log B_\lambda
 \quad\text{and}\quad B_\lambda^{10}\le e^{t/10}.
\]
By Lemma~\ref{lem:easy_algebraic-input} we have
$t\le B_\lambda$ and $B_\lambda\ge3$. Thus $t>100$, hence $B_\lambda>100$ and
$t>100\log100>400$. Using $k\le20(1+t)$ and
$\chi_\lambda\le2(1-\lambda)=2e^{-t}$, we obtain
\begin{align*}
 \chi_\theta B_\theta^{10}
 &\le k^{11}\chi_\lambda B_\lambda^{10} \\
 &\le 2\cdot20^{11}(1+t)^{11}e^{-9t/10} \\
 &\le 2\cdot20^{11}\cdot401^{11}e^{-360}<10^{-100},
\end{align*}
where we used that $(1+t)^{11}e^{-9t/10}$ is decreasing for $t\ge400$.
This concludes the proof.
\end{proof}

\subsection{The reduction to bounded approximation}\label{sec:detail}

We now show how Lemma~\ref{lem:weak} is used. Indeed, 
Lemma~\ref{lem:weak} reduces Proposition~\ref{prop:small-gain} to finding a density function 
satisfying the property of the following lemma.

\begin{lemma}\label{lem:bounded-approximation}
    Let $\nu_{\lambda}$ be absolutely continuous with density function $f_{\lambda}$. Assume that there is a function $g_{\lambda}$ satisfying \begin{equation}\label{eq:analytic-target}
 \|f_\lambda-g_\lambda\|_1\le\frac1{2\cdot21^2}\quad \text{ and } \quad
 \|g_\lambda\|_\infty\le\tfrac12.
\end{equation}
Then
\begin{equation}
 f_\lambda\in L^{1+\delta_\lambda}(\R)\quad \text{ for } \quad
 \delta_\lambda=\frac{1}{10(1-\log(1-\lambda))}.
\end{equation}
\end{lemma}

\begin{proof}
Normalize the support of $\nu_{\lambda}$ to $I=[0,1]$ and denote the resulting
measure and functions by $\sigma=F_{\lambda}\,dx$ and $G_{\lambda}$.
Since the original support has length $2/(1-\lambda)$, we have
\[
 \|F_{\lambda}-G_{\lambda}\|_1\le\frac1{2\cdot21^2}\quad \text{ and } \quad
 \|G_{\lambda}\|_\infty\le K:=(1-\lambda)^{-1}.
\]
Thus $\omega(t):=\sup_{|E|\le t}\sigma(E)\le(2\cdot21^2)^{-1}+Kt$.
By \eqref{eq:local}, applied to $\sigma$ with $c=2/\lambda$,
the hypothesis of Lemma~\ref{lem:weak} holds with
\[
 \eta=21^{-2}\quad \text{ and } \quad \varepsilon=\frac{\lambda}{4\cdot21^2K}.
\]
Consequently, $F_{\lambda}\in L^p$ with
\[
 p-1=\frac{\log21}{2\log(4\cdot21^3K/\lambda)}
 \ge\frac1{10(1+\log K)}=\delta_\lambda,
\]
where we used $\lambda>1/2$ and $K\ge1$.
Compact support and affine rescaling give \eqref{eq:small-gain}.
\end{proof}

Next we recall the definitions of detail and higher order detail from \cites{Kit24, Kit25}. 
We will use bounds on detail to establish the 
assumptions of Lemma~\ref{lem:bounded-approximation} and 
then deduce Proposition~\ref{prop:small-gain}.
Let $\eta_y$ be the centred Gaussian density of variance $y$. We will write
$g_\lambda=\nu_\lambda*\eta_1$, whose supremum is at most
$(2\pi)^{-1/2}<1/2$. Write $\eta_y^{(k)}=\partial_y^k\eta_y$ and set
\[
 Q=\|\eta_1^{(1)}\|_1^{-1}=\sqrt{\pi e/2}\in(2,21/10).
\]
For a probability measure $\sigma$ on $\R$, $r>0$ and an integer $k\ge1$,
define
\[
 s_r^{(k)}(\sigma)
 =r^{2k}Q^k\|\sigma*\eta_{kr^2}^{(k)}\|_1
 \quad \text{ and } \quad s_r(\sigma)=s_r^{(1)}(\sigma).
\]
We refer the reader to \cites{Kit24, Kit25} and \cite{KK25} for a discussion of 
detail and its properties. We just highlight that for $0 < a < b$,
\begin{equation}\label{eq:heat-cost}
 \|\sigma*\eta_{a^2}-\sigma*\eta_{b^2}\|_1
 \le\frac2Q\int_a^b s_t(\sigma)\,\frac{dt}{t}.
\end{equation}
We first establish the following inequality,  which is a more efficient
version of \cite{Kit25}*{Lemma~1.14} and of
\cite{KK25}*{Proposition~4.5}.

\begin{lemma}\label{lem:integrating-detail}
    Suppose that $\sigma$ is a probability measure on $\R$, that $k$ is a positive integer and that $r, b > 0$ with $b \geq r / \sqrt{k}$. Then
    \[
        s_r(\sigma)\le10^k\left(
 \sup_{r/\sqrt{k}\le t\le b}s_t^{(k)}(\sigma)+\frac{r^2}{b^2}\right).
    \]
\end{lemma}

\begin{proof}
The case $k=1$ is immediate so assume $k\ge2$. Repeated integration by parts gives
\begin{align*}
 \frac{(-1)^{k-1}}{(k-2)!}
 \int_{r^2}^{\infty}(y-r^2)^{k-2}\sigma*\eta_y^{(k)}\,dy
 &= \frac{(-1)^{k-2}}{(k-3)!}
 \int_{r^2}^{\infty}(y-r^2)^{k-3}\sigma*\eta_y^{(k-1)}\,dy \\
 &= \dots \\
 &= -\int_{r^2}^{\infty}\sigma*\eta_y^{(2)}\,dy \\
 &= \sigma*\eta_{r^2}^{(1)}.
\end{align*}
Since $\frac{d}{dy}(1-\frac{r^2}{y} )^{k-1}
 =(k-1)r^2(y-r^2)^{k-2}y^{-k},$ the kernel $(k-1)r^2(y-r^2)^{k-2}y^{-k}\,dy$ is a probability measure on
$[r^2,\infty)$.
For $b\ge r/\sqrt{k}$, its mass above $kb^2$ is, by Bernoulli's inequality,
\[
 1-\left(1-\frac{r^2}{kb^2}\right)^{k-1}
 \le\frac{k-1}{k}\frac{r^2}{b^2}
 \le\frac{r^2}{b^2}.
\]
Since
$Q^{1-k}k^k/(k-1)!<10^k$, splitting at $kb^2$ and using that $s_r^{(k)}(\sigma)\le1$ yields
\[
 s_r(\sigma)\le10^k\left(
 \sup_{r/\sqrt{k}\le t\le b}s_t^{(k)}(\sigma)+\frac{r^2}{b^2}\right),
\] concluding the proof.
\end{proof}

Towards applying Lemma~\ref{lem:bounded-approximation}, we note that the following estimates 
of detail are sufficient.

\begin{lemma}[Integrating detail]\label{lem:detail-to-density}
Let $R_0\ge100$ and $\delta>0$. Suppose that a probability measure
$\sigma$ satisfies for $R\ge R_0,\ 1\le k\le10\log R$ and $r\in (0,1)$ that
\begin{equation}\label{eq:higher-detail-hypothesis}
 \begin{aligned}
 s_{e^{-R}}^{(k)}(\sigma)\le100^{-k} \quad\text{ and } \quad
 s_r(\sigma)\le\delta.
 \end{aligned}
\end{equation}
Then $\sigma$ is absolutely continuous with density $f$ and 
\begin{equation}\label{eq:bounded-approximation}
 \|f-\sigma*\eta_1\|_1\le4\delta R_0+\frac1{4R_0}.
\end{equation}
\end{lemma}

\begin{proof}
We will let $R\ge4R_0$ and apply Lemma \ref{lem:integrating-detail} with $r=e^{-R}$, $k=\lceil2\log R\rceil$ and
$b=e^{-R/2}$. 
For $t\in[r/\sqrt{k},b]$, we have $\log(1/t)\ge R/2\ge R_0$
and $k\le10\log(R/2)$. Thus \eqref{eq:higher-detail-hypothesis}
applies uniformly on this interval, bounding the supremum by $100^{-k}$.
Therefore, for $R \geq 4R_0$,
\begin{equation}\label{eq:detail-decay}
 s_{e^{-R}}(\sigma)\le10^k(100^{-k}+e^{-R})\le R^{-2}.
\end{equation}
By \cite{Kit24}*{Lemma~1.18}, $\sigma$ is absolutely continuous and denote by $f\in L^1(\R)$
the density of $\sigma$.
Since $\sigma*\eta_{a^2}\to f$ in $L^1$ as $a\to 0$ and letting
$a\to 0$ in \eqref{eq:heat-cost} with $b=1$, gives
\begin{align*}
 \|f-\sigma*\eta_1\|_1
 &\le\frac2Q\int_0^1s_t(\sigma)\,\frac{dt}{t} \\
 &\le\frac2Q\left(\int_0^{4R_0}\delta\,dR
       +\int_{4R_0}^{\infty}\frac{dR}{R^2}\right) \le4\delta R_0+\frac1{4R_0},
\end{align*}
where we substituted $t=e^{-R}$, used $2/Q<1$, the inequality~\eqref{eq:detail-decay} and the assumption. This proves
\eqref{eq:bounded-approximation}.
\end{proof}

In the next subsection, we will establish the estimates necessary for 
Lemma~\ref{lem:detail-to-density} with suitable constants as stated in the following proposition.

\begin{proposition}[Effective detail]\label{prop:effective-detail}
If $\lambda\in(1/2,1)$ is algebraic and $\chi_\lambda B_\lambda^{10}\le10^{-100}$,
then with $R_0=10^6B_\lambda^3$ we have for $R\ge R_0$ and $\ 1\le k\le10\log R$ that
\begin{equation}\label{eq:effective-detail}
 s_{e^{-R}}^{(k)}(\nu_\lambda)\le100^{-k}.
\end{equation}
Moreover, for every $\lambda\in(1/2,1)$ and $r \in (0,1)$,
\begin{equation}\label{eq:uniform-detail}
 s_r(\nu_\lambda)\le40\chi_\lambda^{1/3}.
\end{equation}
\end{proposition}

Having established all the above and in particular assuming Proposition~\ref{prop:effective-detail}, 
we prove Proposition~\ref{prop:small-gain}.

\begin{proof}[Proof of Proposition~\ref{prop:small-gain}]
Set $R_0=10^6B_\lambda^3$. By Lemma~\ref{lem:easy_algebraic-input},
$B_\lambda\ge3$, so $R_0\ge100$.
Proposition~\ref{prop:effective-detail} shows that $\nu_\lambda$ satisfies
the hypotheses of Lemma~\ref{lem:detail-to-density} with this value of $R_0$
and $\delta=40\chi_\lambda^{1/3}$. Thus $\nu_\lambda$ is absolutely continuous
with density $f_\lambda$. Writing $g_\lambda=\nu_\lambda*\eta_1$,
\eqref{eq:bounded-approximation} gives
\begin{align*}
 \|f_\lambda-g_\lambda\|_1
 &\le160\chi_\lambda^{1/3}R_0+\frac1{4R_0} \\
 &\le160\cdot10^6\cdot10^{-100/3}+\frac1{108\cdot10^6}
 <\frac1{2\cdot21^2},
\end{align*}
where we used $B_\lambda\ge3$ from \eqref{eq:arithmetic} and
$\chi_\lambda^{1/3}B_\lambda^3
=(\chi_\lambda B_\lambda^{10})^{1/3}B_\lambda^{-1/3}\le10^{-100/3}$ by the assumption of 
Proposition~\ref{prop:small-gain} that $\chi_\lambda B_\lambda^{10} \leq 10^{-100}$.
Moreover, since $\nu_\lambda$ is a probability measure, $\|g_\lambda\|_\infty\le\|\eta_1\|_\infty
 =(2\pi)^{-1/2}<\tfrac12.$
Thus $f_\lambda$ and $g_\lambda$ satisfy \eqref{eq:analytic-target} and
Lemma~\ref{lem:bounded-approximation} gives
$f_\lambda\in L^{1+\delta_\lambda}(\mathbb R)$, as required.
\end{proof}

\subsection{Effective detail from entropy and variance}\label{sec:effective-detail-proof}

It remains to prove Proposition~\ref{prop:effective-detail}, which is
the goal of this subsection. The proof presented here is a simplified version 
of the variance summation method as developed in \cite{KK25}. To deduce
absolute continuity as in \cite{KK25} it suffices to prove results about 
the measure in question for sufficiently small scales. In contrast, 
Proposition~\ref{prop:effective-detail} requires bounds at the explicit scales
$e^{-R}$ for $R \geq R_0$ and $R_0 = 10^6 B_{\lambda}^3$. The latter is necessary in order to apply 
Lemma~\ref{lem:bounded-approximation}, which relies on Lemma~\ref{lem:weak}. 
To prove Proposition~\ref{prop:effective-detail} we therefore require explicit
separation estimates, which we prove in Lemma~\ref{lem:algebraic-input} and compare
to previous separation estimates right after stating Lemma~\ref{lem:algebraic-input}.

We first bound detail
in terms of conditional variance as stated in Lemma~\ref{lem:conditional-variance}. 
Then, we use entropy and algebraic separation
to find enough conditional variance in disjoint blocks of the random power series 
determining
the Bernoulli convolution.

We will use for probability measures $\sigma_1, \ldots , \sigma_k$ that, 
as shown in \cite{KK25}*{Lemma~4.3},
\begin{equation}\label{eq:products}
 \begin{aligned}
 s_r^{(k)}(\sigma_1*\cdots*\sigma_k)
 &\le\prod_{i=1}^k s_r(\sigma_i)
 \quad \text{ and } \quad s_r(\mathcal N(0,v))=\frac{r^2}{r^2+v},
 \end{aligned}
\end{equation} where $\mathcal N(0,v)$ is the standard Gaussian measure and $v \geq 0$.

First, we establish Lemma~\ref{lem:small-summands},
which gives an explicit one-dimensional version of
\cite{KK25}*{Lemma~4.11}. We write $\mathcal W_1$ for the $1$-Wasserstein
distance, applied to random variables through their laws and use the following
Berry--Esseen estimate.

\begin{lemma}[Berry--Esseen estimate, \cite{Gol10}*{Theorem~1.1}] \label{lem:berry_essen}
    Let $\varepsilon > 0$ and let $U_1,\ldots,U_n$ be independent real
random variables satisfying $|U_i|\le\varepsilon$. Let $\mu_i = \mathbb{E}[U_i]$ and let $\omega_i^2 = \Var[U_i]$.
 Let $\mu = \sum_{i=1}^n \mu_i$ and let $\omega^2 = \sum_{i=1}^n \omega_i^2$. 
 Then $$\mathcal{W}_1\left(\sum_{i=1}^n U_i, N(\mu, \omega^2)\right) < 2 \varepsilon.$$
\end{lemma}

\begin{lemma}[Small summands]\label{lem:small-summands}
Let $r>0$, $\varepsilon>0$, and let $U_1,\ldots,U_n$ be independent real
random variables satisfying $|U_i|\le\varepsilon r$ and write
$v=\sum_i\Var(U_i)$. Then
\[
 s_r\left(\sum_iU_i\right)
 \le \frac{r^2}{r^2+v}+6\varepsilon.
\]
\end{lemma}

\begin{proof}
    By \cite{KK25}*{Lemma~4.7}, for probability measures $\lambda_1,\lambda_2$ on $\R$,
    \[
        |s_r(\lambda_1)-s_r(\lambda_2)|
        \le e r^{-1}\mathcal W_1(\lambda_1,\lambda_2).
    \]
    Let $\mu = \sum_i \mathbb{E}[U_i]$ and let $\lambda$ be the law of $N(\mu, v)$.
    Together with Lemma~\ref{lem:berry_essen}, this gives
    \begin{align*}
        s_r \left(\sum_iU_i\right) & \leq s_r\left(\lambda\right) + \left| s_r \left(\sum_iU_i\right) - s_r\left(\lambda\right) \right| \\
        & \leq \frac{r^2}{r^2+v} + e r^{-1} \mathcal{W}_1(\lambda, \sum_iU_i) \\
        & \leq \frac{r^2}{r^2+v} + 6 \varepsilon.
    \end{align*}
    This concludes the proof.
\end{proof}

We can already prove \eqref{eq:uniform-detail}. Let $\lambda\in(1/2,1)$,
$r\in(0,1)$ and $\varepsilon=\chi_\lambda^{1/3}<1$. Choose $n\ge0$ with
$\lambda\varepsilon r<\lambda^n\le\varepsilon r$ and define
$$Y_N=\sum_{j=n}^{N}\xi_j\lambda^j$$ for $\xi_j$ independent random variables
satisfying $\mathbb{P}(\xi_j = 1) = \mathbb{P}(\xi_j = -1) = \frac{1}{2}$. Denote
\[
 v_N:=\Var(Y_N)=\frac{\lambda^{2n}-\lambda^{2N+2}}{1-\lambda^2}
\] and since
$1-\lambda^2\le2\chi_\lambda$ and $\lambda^2>1/4$, $$\lim_{N\to \infty} v_N = \frac{\lambda^{2n}}{1-\lambda^2} 
 > \frac{\varepsilon^2r^2}{8\chi_\lambda}
 =\frac{r^2}{8\chi_\lambda^{1/3}}.$$
Applying Lemma~\ref{lem:small-summands} for all sufficiently large $N$ we have
\begin{equation}\label{eq:uniform-detail-proof}
 s_r(\nu_\lambda)\le s_r(Y_N)
 \le\frac{r^2}{r^2+v_N}+6\chi_\lambda^{1/3}
 <40\chi_\lambda^{1/3}.
\end{equation}

We next prove the following version of \cite{KK25}*{Proposition~4.8}, which 
is an easy consequence of \eqref{eq:products} and Lemma \ref{lem:small-summands}.

\begin{lemma}[Conditional variance]\label{lem:conditional-variance}
Denote $\varepsilon_0=10^{-4}$, let $r>0$ and let $k\ge1$ be an integer.
Let $X_1,\ldots,X_n$ be real random variables with finite second moments
on a probability space $(\Omega,\mathcal F,\Prob)$, and let
$\mathcal A_1,\ldots,\mathcal A_n$ be sub-$\sigma$-algebras of $\mathcal F$.
Suppose that the $\sigma$-algebras $\sigma(X_i)\vee\mathcal A_i$ are
independent. Assume that there are $\mathcal A_i$-measurable real random
variables $C_i$ such that
\[
 |X_i-C_i|\le\varepsilon_0r
 \quad \text{ and } \quad
 \sum_i\E[\Var(X_i\mid\mathcal A_i)]\ge1000kr^2.
\]
Then
\[
 s_r^{(k)}\left(\sum_iX_i\right)\le100^{-k}.
\]
\end{lemma}

\begin{proof}
Write $W_i=\Var(X_i\mid\mathcal A_i)$. Note that the $W_i$ are $\mathcal{A}_i$-measurable random variables. 
The $W_i$ are independent and
$0\le W_i\le\E[(X_i-C_i)^2\mid\mathcal A_i]\le\varepsilon_0^2r^2$ almost surely. In particular,

\begin{equation}
    \E [W_i] \in [0, \varepsilon_0^2r^2] \quad \text{ and } \quad \Var [W_i] \leq \mathbb{E}[W_i^2] \leq \varepsilon_0^2r^2 \E [W_i]. \label{eq:w_bounded}
\end{equation}

By \eqref{eq:w_bounded} and a greedy algorithm we can partition $\{1, \dots, n\}$ into $k$ sets $J_1, \dots, J_k$ such that
$$
\sum_{j \in J_i} \E [W_j] \geq 999 r^2
$$

For one such $J_i$, write $m_i=\sum_{j\in J_i}\E[W_j]\ge999r^2$.
Since by \eqref{eq:w_bounded} we have $\Var(\sum_{j\in J_i}W_j)\le\varepsilon_0^2r^2m_i$, Chebyshev's
inequality gives
\[
 \Prob\left(\sum_{j\in J_i}W_j<m_i/2\right)
 \le\frac{4\varepsilon_0^2r^2}{m_i}\le\frac{4\varepsilon_0^2}{999}.
\]
Condition on $\bigvee_{j\in J_i}\mathcal A_j$. The $X_j$ remain
conditionally independent, with conditional variances $W_j$, and the
centres $C_j$ are measurable. On the event $\sum_{j\in J_i}W_j\ge m_i/2$,
Lemma~\ref{lem:small-summands} bounds the conditional detail by
$2/999+6\varepsilon_0$ and on its complement use the bound $1$.
Convexity and translation invariance therefore give
\[
 s_r\left(\sum_{j\in J_i}X_j\right)
 \le\frac2{999}+6\varepsilon_0+\frac{4\varepsilon_0^2}{999}<\frac1{100}.
\]
The $k$ group sums are independent. Apply \eqref{eq:products} to obtain the claimed bound.
\end{proof}

We use Varj\'u's averaged entropy at scale \cite{Var19}*{Section~2.2},
though unlike Varj\'u we use natural logarithms instead of base 2. For a finitely supported real random variable
$X$ and $s>0$, define
\[
 H(X;s)=\frac1s\int_0^sH\left(\left\lfloor \frac{X+t}{s}\right\rfloor\right)\,dt,
\]
where $H(\cdot)$ is the Shannon entropy.
If $T$ is independent of $X$ and uniform on $[0,1]$, write
\begin{equation}\label{eq:variance-at-scale}
 \mathcal A_s=\sigma(T,\lfloor \tfrac{X}{s}+T\rfloor) \quad \text{ and } \quad
 \qquad
 V(X;s)=s^{-2}\E[\Var(X\mid\mathcal A_s)].
\end{equation}
The $\mathcal A_s$-measurable midpoint
$C_s=s(\lfloor \frac{X}{s} +T\rfloor-T+1/2)$ satisfies $|X-C_s|\le s/2$.
We note that $V(X;s)$ is an additive analogue, with a fixed conditioning, of the
trace in \cite{KK25}*{Section~5.1}.

As in \cite{KK26}*{Theorem~1.4} and as applied in
\cite{KK25}*{(5.2)}, an entropy gap forces conditional
variance. In this paper we use the following elementary bound.

\begin{lemma}[Entropy gives variance]\label{lem:entropy-variance}
For every finitely supported $X$ and $s>0$,
\[
 H(X;s/4)-H(X;s)\le10V(X;s)^{1/4}.
\]
\end{lemma}

We compare this lemma to \cite{KK26}*{Theorem~1.4}. 
There, the authors show that, with their notion of entropy at a scale $H_a(\cdot ; s)$ 
which behaves in a qualitatively similar way to Varj\'u's entropy at a scale, we have
for the trace $\mathrm{tr}(X;s)$ as defined in \cite{KK26}, 
which is analogous to $V(X;s)$, that
$$
\mathrm{tr}(X;s) \gg a^{-2} (H_a(X;s/4a) - H_a(X; s / a) - O(e^{-a^2 / 4} + a^3 s)).
$$
By substituting $a = C \sqrt{ \log (H_a(X;s/4a) - H_a(X; s / a))^{-1}}$ 
we deduce
$$
\mathrm{tr}(X; s) \gg \frac{\tilde{\alpha}}{\log \tilde{\alpha}^{-1}}
$$
where $\tilde{\alpha} = H_a(X;s/4a) - H_a(X; s / a)$ which, at least when the gap is small, is stronger than the estimate of
$V(X; s) \geq \alpha^4 / 10^4$ for $\alpha = H(X;s/4)-H(X;s)$ coming from Lemma \ref{lem:entropy-variance}. The estimate from Lemma \ref{lem:entropy-variance} is still
useful for our current purposes as it is simple and gives reasonable bounds at all scales.

\begin{proof}[Proof of Lemma \ref{lem:entropy-variance}]
First, let $Y$ be finitely supported in an interval of length $s$. A shifted grid
of mesh $s/4$ has at most five occupied cells. For a probability vector
$p=(p_1,\ldots,p_m)$ with $m\le5$, the bound
$-\log u\le2(u^{-1/2}-1)$ implies
$-p_j\log p_j\le2(\sqrt{p_j}-p_j)\le2\sqrt{p_j}(1-p_j)$.
Since $\sum_j(1-p_j)=m-1$, Cauchy--Schwarz applied to the factors
$\sqrt{p_j(1-p_j)}$ and $\sqrt{1-p_j}$ therefore gives
\[
 H(p)\le2\sum_j\sqrt{p_j}(1-p_j)
 \le2\sqrt{m-1}\sqrt{1-\sum_jp_j^2}
 \le4\sqrt{1-\sum_jp_j^2}.
\]
Let $Y'$ be an independent copy of $Y$. For a fixed grid,
$1-\sum_jp_j^2$ is the probability that $Y$ and $Y'$ occupy different cells.
Averaging over a uniform shift, this probability becomes
$\E[\min\{1,4|Y-Y'|/s\}]$. Jensen's and Cauchy--Schwarz inequalities yield
\[
 H(Y;s/4)
 \le4\left(\frac4s\E[|Y-Y'|]\right)^{1/2}
 \le8\cdot2^{1/4}\left(\frac{\Var(Y)}{s^2}\right)^{1/4}
 \le10\left(\frac{\Var(Y)}{s^2}\right)^{1/4}.
\]
Now take uniform shifts $T,U$, independent of each other and of $X$,
for the coarse and fine grids.
Writing $A=\lfloor X/s+T\rfloor$ and $Z=\lfloor4X/s+U\rfloor$, the chain
rule gives
\[
 H(X;s/4)-H(X;s)
 =H(Z\mid T,U)-H(A\mid T,U)
 \le H(Z\mid A,T,U).
\]
Conditionally on $(A,T)$, the variable $X$ lies in an interval of length
$s$ and $U$ remains independent. Applying the preceding bound to this
conditional law and then Jensen's inequality gives
\[
 H(Z\mid A,T,U)
 \le10\E\left[\left(\frac{\Var(X\mid\mathcal A_s)}{s^2}\right)^{1/4}\right]
 \le10V(X;s)^{1/4}.
\]
This concludes the proof.
\end{proof}

We next record the arithmetic estimates needed to produce an entropy gap. Recall that we defined 
$$B_{\lambda} = \frac{1 + d_{\lambda} + \log M_{\lambda}}{\min\{ 1, \log M_{\lambda} \}}.$$ 
From now on we also denote for an integer $\ell\ge1$
by $X_{\ell}$ the random variable $$X_\ell=\sum_{j=0}^{\ell-1}\xi_j\lambda^j,$$ 
where the $\xi_j$
are independent 
and satisfy $\mathbb{P}(\xi_j = 1) = \mathbb{P}(\xi_j = -1) = \frac{1}{2}$. 

\begin{lemma}[Arithmetic input]\label{lem:algebraic-input}
Denote by $\Delta_\ell$ the minimum distance between distinct atoms of $X_{\ell}$. Then
\begin{equation}\label{eq:arithmetic-estimates}
 \Delta_\ell\ge (\ell+1)^{1-d_{\lambda}}M_\lambda^{-\ell} \ge e^{-B_\lambda\ell}\quad \text{ and } \quad
 H(X_\ell)\ge\frac{\ell\log2}{B_\lambda}.
\end{equation}
\end{lemma}

We discuss how this lemma compares to previous separation estimates. 
Indeed, Garsia's separation estimate \cite{Gar62}*{Lemma~1.51} gives
$\Delta_\ell\ge c_\lambda\ell^{-k_\lambda}M_\lambda^{-\ell},$
where $k_\lambda$ is the number of conjugates of $\lambda$ on the unit circle
and $c_\lambda>0$ is a fixed constant depending only on $\lambda$.
Unlike in \cite{KK25}, \cite{Var19}, \cite{Kit24}, \cite{Kit25},
in this paper we care about the strength of separation at small scales 
instead of just the asymptotic rate of exponential separation. 
Thus, to use Garsia's estimate 
we would need to bound $c_{\lambda}$. 
Our elementary separation bound avoids this additional estimate and gives an 
explicit bound for every $\ell$ directly in terms of \(d_\lambda\) and \(M_\lambda\). 
The estimate of $\Delta_\ell$ in Lemma \ref{lem:algebraic-input} 
has the same asymptotic rate of exponential separation as Garsia's estimate 
although the exact separation predicted is slightly weaker for large $\ell$.

For the entropy $H(X_{\ell})$, Breuillard--Varj\'u \cite{BV20}*{Theorem~5} give the stronger,
degree-independent estimate
\begin{equation}\label{eq:bv-entropy}
 H(X_\ell)\ge 0.3\cdot\ell\min\{1,\log M_\lambda\},
\end{equation}
where the constant $0.3$ in \eqref{eq:bv-entropy} is based on the non-rigorous
numerical estimate of \cite{BV20}*{Theorem~5}. The logarithms in \cite{BV20} are base $2$, so the bound 
$0.44$ as stated in \cite{BV20}*{Theorem~5} gives in natural logarithms $0.44\log2\ge0.3$.
The entropy estimate in Lemma~\ref{lem:algebraic-input} is weaker than
\eqref{eq:bv-entropy} due to the dependence on the degree $d_{\lambda}$,
yet has a simple proof. Nonetheless, using \eqref{eq:bv-entropy}
would not remove the degree dependence
in Theorem~\ref{thm:algebraic} as we need to understand the entropy
at small explicit scales which requires suitable explicit bounds on $\Delta_{\ell}$.

\begin{proof}[Proof of Lemma~\ref{lem:algebraic-input}]
Let $P_\lambda(T)=a\prod_{i=1}^{d_\lambda}(T-\lambda_i)$ be the primitive
integer minimal polynomial, with $a>0$ and $\lambda_1=\lambda$.
For the separation bound in \eqref{eq:arithmetic-estimates}, suppose that $Q$ is a degree $n<\ell$
polynomial and
has coefficients in $\{-1,0,1\}$.  Also assume that $Q(\lambda)\ne0$.
Then the resultant of $P_\lambda$ and $Q$ is a nonzero integer, so 
$1\le a^n\prod_{i=1}^{d_\lambda}|Q(\lambda_i)|.$
Using $|Q(\lambda_i)|\le(n+1)\max\{1,|\lambda_i|\}^n$ for $i\ge2$ gives
$|Q(\lambda)|\ge(n+1)^{1-d_\lambda}M_\lambda^{-n}.$
Differences between distinct atoms of $X_\ell$ are twice such values, hence
\[
 \Delta_\ell\ge2\ell^{1-d_\lambda}M_\lambda^{-(\ell-1)}
 \ge(\ell+1)^{1-d_\lambda}M_\lambda^{-\ell}.
\]
Since $\log(\ell+1)\le\ell$ and
$B_\lambda\ge d_\lambda-1+\log M_\lambda$, the last expression is at least
$e^{-B_\lambda\ell}$.

We now prove the entropy bound in \eqref{eq:arithmetic-estimates}.
If $\lambda$ is not an algebraic unit, it cannot be a root of a nonzero
polynomial with coefficients in $\{-1,0,1\}$. Indeed,
after factoring out a power,
such a polynomial has leading and constant coefficients in $\{-1,1\}$
and therefore any zero of such a polynomial is an algebraic unit.
Thus, each choice of signs gives a different partial sum and so
$H(X_\ell)=\ell\log2\ge\ell\log2/B_\lambda$.

We may therefore assume that $\lambda$ is an algebraic unit. In this case,
$\log M_\lambda
 =\sum_{i=1}^{d_\lambda}\log\max\{1,|\lambda_i|\},$
so there is a conjugate $z$ of $\lambda$ with
$\log|z|\ge\log M_\lambda/d_\lambda$.
Let $k=\lceil d_\lambda\log2/\log M_\lambda\rceil$.
Then $|z|^k\ge2$ and, by the definition of $B_\lambda$,
\[
 k\le1+\frac{d_\lambda\log2}{\log M_\lambda}\le B_\lambda.
\]
Denote $m=\lceil\ell/k\rceil$ and
$Y=\sum_{j=0}^{m-1}\xi_{jk}\lambda^{jk}$.
We claim that the $2^m$ choices of signs give distinct values of $Y$.
Indeed, an equality between two such sums would also hold with $\lambda$
replaced by $z$ as any polynomial $Q$ defining such an equality 
would satisfy that $P_{\lambda}$ divides $Q$. 
As $|z|^k \geq 2$ it holds for all integers $t\geq 1$ that
$|z|^{kt}>\sum_{j=0}^{t-1}|z|^{kj}$ and therefore different signs always lead to different values.
The remaining digits are independent and therefore increase entropy, showing that
\[
 H(X_\ell)\ge H(Y)=m\log2
 \ge\frac{\ell\log2}{B_\lambda},
\]
as required.
\end{proof}

We now turn the entropy gap into variance over a range of scales,
which is analogous to \cite{KK25}*{Propositions~5.1 and~5.4}.

\begin{lemma}[Variance in a block]\label{lem:block-variance}
Let $X_{\ell}$ and $B_{\lambda}$ be defined as above. Then 
every integer $\ell\ge100B_\lambda^2$ satisfies
\begin{equation}\label{eq:block-variance}
 \int_{e^{-2B_\lambda\ell}}^{e^{-\ell/(32B_\lambda)}}V(X_\ell;s)\,\frac{ds}{s}
 \ge\frac{\ell}{10^{10}B_\lambda^7}.
\end{equation}
\end{lemma}

\begin{proof}
Write $a=e^{-2B_\lambda\ell}$, $b=e^{-\ell/(32B_\lambda)}$ and denote
$$ D=\int_a^b\bigl(H(X_\ell;s/4)-H(X_\ell;s)\bigr)\frac{ds}{s}.$$
Since $\log(b/a)\le2B_\lambda\ell$, Lemma~\ref{lem:entropy-variance}
and H\"older's inequality with exponents $4$ and $4/3$ give
\begin{align}
 D\le10\int_a^b V(X_\ell;s)^{1/4}\,\frac{ds}{s} &\leq10\left( \int_a^b 1 \frac{ds}{s} \right)^{3/4} \left(\int_a^b V(X_\ell;s)\,\frac{ds}{s}\right)^{1/4}\notag
 \\ &\le10(2B_\lambda\ell)^{3/4}
 \left(\int_a^b V(X_\ell;s)\,\frac{ds}{s}\right)^{1/4}.\label{eq:block-variance-holder}
\end{align}
To estimate $D$ from below note that 
$$D=\int_{a/4}^aH(X_\ell;t)\frac{dt}{t} -\int_{b/4}^bH(X_\ell;t)\frac{dt}{t}.$$ So we need
to estimate $H(X_\ell;t)$ in the ranges $t \in [a/4,a]$ and $t \in [b/4,b]$.
By \eqref{eq:arithmetic-estimates}, $\Delta_\ell\ge e^{-B_\lambda\ell}>a$,
so $H(X_\ell;t)=H(X_\ell)\ge\ell\log2/B_\lambda$ for $t\le a$. It remains to consider
the range $t \in [b/4,b]$.
Since $|X_\ell|\le(1-\lambda)^{-1}$, a grid of mesh $t$ meets the support of $X_{\ell}$
in at most $2(1-\lambda)^{-1}t^{-1}+2$ cells.
For $t \in [b/4,b]$ we therefore have
\begin{align*}
 H(X_\ell;t)
 \le\log\bigl(2(1-\lambda)^{-1}t^{-1}+2\bigr) &\leq \log\bigl(10(1-\lambda)^{-1}b^{-1}\bigr) \\
 &\le B_\lambda+\frac{\ell}{32B_\lambda}+\log10\le\frac{\ell}{10B_\lambda},
\end{align*}
where we used in the penultimate inequality that $\log (1- \lambda)^{-1} \leq B_{\lambda}$ 
from Lemma~\ref{lem:easy_algebraic-input}  and in the last inequality 
the assumption $\ell\ge100B_\lambda^2$ and that $B_\lambda\ge3$ 
from Lemma~\ref{lem:easy_algebraic-input}.
Thus,
\begin{align*}
 D=\int_{a/4}^aH(X_\ell;t)\frac{dt}{t}
  -\int_{b/4}^bH(X_\ell;t)\frac{dt}{t}
  \ge\frac{(\log4)(\log2-1/10)\ell}{B_\lambda} \geq \frac{\ell}{10 B_{\lambda}}.
\end{align*}
Applying the latter inequality for $D$ 
in \eqref{eq:block-variance-holder} and rearranging concludes the proof as
\[
 \int_a^b V(X_\ell;s)\,\frac{ds}{s}
 \geq \frac{(\ell/10B_\lambda)^4}{10^4(2B_\lambda\ell)^3}
 \geq \frac{\ell}{8 \cdot 10^{8}B_\lambda^7} \geq \frac{\ell}{10^{10}B_\lambda^7} .\qedhere
\]
\end{proof}

We finally prove Proposition~\ref{prop:effective-detail},
which resembles \cite{KK25}*{Proposition~6.7}.

\begin{proof}[Proof of Proposition~\ref{prop:effective-detail}]
We have already proved \eqref{eq:uniform-detail} in \eqref{eq:uniform-detail-proof}. For
\eqref{eq:effective-detail}, fix $r=e^{-R}$ with $R\ge R_0=10^6B_\lambda^3$. Denote
$\varepsilon_0=10^{-4}$ and $N=\lceil1000B_\lambda^2\rceil$.
We will apply the previous results for $X_{\ell}$ with $\ell=N,N^2,\ldots$ 
with, as will become apparent later, $\ell \leq R/(4B_\lambda)$.
Each length satisfies $\ell\ge100B_\lambda^2$, as required in
Lemma~\ref{lem:block-variance}.

The goal is to apply Lemma~\ref{lem:conditional-variance}. 
First, fix one choice of $\ell$ and set
$I_\ell=[\ell/(32B_\lambda),2B_\lambda\ell]$ and $h=\ell\chi_\lambda$.
Note that
\begin{equation}\label{eq:block-variance-grid}
 \frac1h\int_0^h
 \sum_{v\in(u+h\mathbb Z)\cap I_\ell}V(X_\ell;e^{-v})\,du
 =\frac1h\int_{I_\ell}V(X_\ell;e^{-v})\,dv
 \ge\frac1{10^{10}\chi_\lambda B_\lambda^7},
\end{equation}
by Lemma~\ref{lem:block-variance}. Choose a shift $u_\ell \in [0,h]$ for which
the sum in \eqref{eq:block-variance-grid} is at least its average, and write
$\mathcal V_\ell=(u_\ell+h\mathbb Z)\cap I_\ell$.
For each $v\in\mathcal V_\ell$, set $s(v)=e^{-v}$ and denote
\[
 n(v)=\left\lceil\frac{R+\log\varepsilon_0^{-1}-v}{\chi_\lambda}\right\rceil.
\]
As $v\le2B_\lambda\ell\le R/2$, we have $n(v)\ge0$ and observe
\begin{equation}\label{eq:block-placement}
 \lambda\varepsilon_0 r\le\lambda^{n(v)}s(v)\le\varepsilon_0 r.
\end{equation}
Use the digits $\xi_j$ starting at $n(v)$ to define
\[
 X_{\ell,v}=\sum_{j=0}^{\ell-1}\xi_{n(v)+j}\lambda^j \quad \text{ and } \quad
Y_{\ell,v}=\lambda^{n(v)}X_{\ell,v}.
\]
The variable $X_{\ell,v}$ has the same law as $X_\ell$.
Consecutive elements of $\mathcal V_\ell$ differ by $\ell\chi_\lambda$, so
their starting indices differ by exactly $\ell$. Therefore the blocks of
this length use disjoint digits $\xi_j$ and are independent.

For each block, take a random shift $T_{\ell,v}$ that is uniform on $[0,1]$, with all shifts
independent of one another and of the digits. Then we condition on
\[
 \mathcal A_{\ell,v}
 =\sigma(T_{\ell,v},\lfloor X_{\ell,v}/s(v)+T_{\ell,v}\rfloor).
\]
As explained after \eqref{eq:variance-at-scale}, the midpoint of the scaled interval gives an
$\mathcal A_{\ell,v}$-measurable centre
\[
 C_{\ell,v}=\lambda^{n(v)}s(v)
 \left(\left\lfloor\frac{X_{\ell,v}}{s(v)}+T_{\ell,v}\right\rfloor
 -T_{\ell,v}+\frac12\right)
\]
with $|Y_{\ell,v}-C_{\ell,v}|\le\lambda^{n(v)}s(v)/2\le\varepsilon_0 r$.
Moreover, the definition of $V(X_\ell;e^{-v})$ and \eqref{eq:block-placement} give
\[
 \begin{aligned}
 \sum_{v\in\mathcal V_\ell}
 \E[\Var(Y_{\ell,v}\mid\mathcal A_{\ell,v})]
 &=\sum_{v\in\mathcal V_\ell}\lambda^{2n(v)}e^{-2v}V(X_\ell;e^{-v})\\
 &\ge\frac{\varepsilon_0^2r^2}{4\cdot10^{10}\chi_\lambda B_\lambda^7},
 \end{aligned}
\]
where we used \eqref{eq:block-variance-grid} and $\lambda^2>1/4$.

We next check that the blocks $Y_{\ell,v}$ for different $\ell = N, N^2\ldots$ 
use disjoint digits $\xi_j$ of $\sum_{j= 0}^{\infty} \xi_j \lambda^j$.
A digit in $Y_{\ell,v}$ has index $j=n(v)+i$ with $0\le i<\ell$.
By the definition of $n(v)$,
\[
 v-(i+1)\chi_\lambda<R+\log\varepsilon_0^{-1}-j\chi_\lambda
 \le v-i\chi_\lambda.
\]
Since $v\in I_\ell$ and the assumption $\chi_\lambda B_\lambda^{10}\le10^{-100}$
implies $\chi_\lambda\le1/(64B_\lambda)$, it follows that
\[
 \frac{\ell}{64B_\lambda}
 \le\frac{\ell}{32B_\lambda}-\ell\chi_\lambda
 <R+\log\varepsilon_0^{-1}-j\chi_\lambda
 \le2B_\lambda\ell.
\]
For distinct admissible lengths $\ell<\ell'$, we have $\ell'\ge N\ell$.
As $N>128B_\lambda^2$,
\[
 2B_\lambda\ell<\frac{N\ell}{64B_\lambda}
 \le\frac{\ell'}{64B_\lambda}.
\]
Thus the ranges above are disjoint, so no digit index can occur in blocks
of different lengths. Together with disjointness for each fixed length
and independence of the digits and shifts, this shows that the
$\sigma$-algebras $\sigma(Y_{\ell,v})\vee\mathcal A_{\ell,v}$ are independent.

Finally, Lemma~\ref{lem:easy_algebraic-input} gives $B_\lambda\ge3$. 
Since $R\ge10^6B_\lambda^3$, 
it follows that
$R/(4B_\lambda)\ge N$ and
$\log(R/(4B_\lambda))\ge\frac12\log R$. 
Using the coarse bound $\lfloor x\rfloor\ge x/2$
for $x\ge1$, the number of admissible $\ell$ is therefore
\[
 \left\lfloor\frac{\log(R/(4B_\lambda))}{\log N}\right\rfloor
 \ge\frac{\log R}{4\log N}.
\]
Since $N\le1001B_\lambda^2$ and $B_\lambda\ge3$, we also have $\log N\le B_\lambda^3$.
Summing the expected conditional variances of the random variables $Y_{\ell,v}$ yields
\[
 \begin{aligned}
 \sum_{\ell}\sum_{v\in\mathcal V_\ell}
 \E[\Var(Y_{\ell,v}\mid\mathcal A_{\ell,v})]
 &\ge\frac{\varepsilon_0^2\log R}{16\cdot10^{10}\chi_\lambda B_\lambda^7\log N}r^2\\
 &\ge\frac{10^{82}}{16}(\log R)r^2>10^5(\log R)r^2,
 \end{aligned}
\]
using $\varepsilon_0=10^{-4}$ and $\chi_\lambda B_\lambda^{10}\le10^{-100}$.
This is larger than $1000kr^2$ whenever $1\le k\le10\log R$.
Lemma~\ref{lem:conditional-variance} therefore bounds the order $k$ detail
of the sum of the blocks $Y_{\ell,v}$ by $100^{-k}$.
The unused digits $\xi_i$ are independent of these random variables, so convolution with
their sum proves \eqref{eq:effective-detail}.
\end{proof}


\begin{bibdiv}
\renewcommand{\PrintDOI}[1]{DOI \href{https://doi.org/#1}{\nolinkurl{#1}}}
\begin{biblist}

\bib{AHT17}{article}{
 author={Anderson, Theresa C.},
 author={Hyt\"onen, Tuomas},
 author={Tapiola, Olli},
 title={Weak $A_\infty$ weights and weak reverse H\"older property in a space of homogeneous type},
 journal={J. Geom. Anal.},
 volume={27},
 number={1},
 date={2017},
 pages={95--119},
 doi={10.1007/s12220-015-9675-6},
}

\bib{BSS23}{book}{
 author={B\'ar\'any, Bal\'azs},
 author={Simon, K\'aroly},
 author={Solomyak, Boris},
 title={Self-similar and self-affine sets and measures},
 series={Mathematical Surveys and Monographs},
 volume={276},
 publisher={American Mathematical Society},
 address={Providence, RI},
 date={2023},
 doi={10.1090/surv/276},
 note={\PrintDOI{10.1090/surv/276}},
}

\bib{BG06}{book}{
 author={Bombieri, Enrico},
 author={Gubler, Walter},
 title={Heights in Diophantine geometry},
 series={New Mathematical Monographs},
 volume={4},
 publisher={Cambridge University Press},
 address={Cambridge},
 date={2006},
 doi={10.1017/CBO9780511542879},
}

\bib{BV20}{article}{
 author={Breuillard, Emmanuel},
 author={Varj\'u, P\'eter P.},
 title={Entropy of Bernoulli convolutions and uniform exponential growth for linear groups},
 journal={J. Anal. Math.},
 volume={140},
 number={2},
 date={2020},
 pages={443--481},
 doi={10.1007/s11854-020-0100-0},
 eprint={arXiv:1510.04043},
}

\bib{DLS14}{article}{
 author={De Lellis, Camillo},
 author={Spadaro, Emanuele},
 title={Regularity of area minimizing currents I: gradient $L^p$ estimates},
 journal={Geom. Funct. Anal.},
 volume={24},
 date={2014},
 number={6},
 pages={1831--1884},
 doi={10.1007/s00039-014-0306-3},
}

\bib{DX10}{article}{
 author={Dubickas, Art\=uras},
 author={Xu, Zhiqiang},
 title={Refinement equations and spline functions},
 journal={Adv. Comput. Math.},
 volume={32},
 date={2010},
 number={1},
 pages={1--23},
 doi={10.1007/s10444-008-9079-2},
}

\bib{Fal14}{book}{
 author={Falconer, Kenneth},
 title={Fractal geometry: Mathematical foundations and applications},
 edition={3},
 publisher={John Wiley \& Sons},
 date={2014},
}

\bib{Gar62}{article}{
 author={Garsia, Adriano M.},
 title={Arithmetic properties of Bernoulli convolutions},
 journal={Trans. Amer. Math. Soc.},
 volume={102},
 number={3},
 date={1962},
 pages={409--432},
 doi={10.2307/1993615},
}

\bib{GM86}{article}{
 author={Giaquinta, Mariano},
 author={Modica, Giuseppe},
 title={Partial regularity of minimizers of quasiconvex integrals},
 journal={Ann. Inst. H. Poincar\'e Anal. Non Lin\'eaire},
 volume={3},
 date={1986},
 number={3},
 pages={185--208},
 url={https://www.numdam.org/item/AIHPC_1986__3_3_185_0/},
}

\bib{Gol10}{article}{
 author={Goldstein, Larry},
 title={Bounds on the constant in the mean central limit theorem},
 journal={Ann. Probab.},
 volume={38},
 date={2010},
 number={4},
 pages={1672--1689},
 doi={10.1214/10-AOP527},
 eprint={arXiv:0906.5145},
}

\bib{Gra14}{book}{
 author={Grafakos, Loukas},
 title={Classical Fourier analysis},
 edition={3},
 series={Graduate Texts in Mathematics},
 volume={249},
 publisher={Springer},
 address={New York},
 date={2014},
 doi={10.1007/978-1-4939-1194-3},
 note={\PrintDOI{10.1007/978-1-4939-1194-3}},
}

\bib{Hut81}{article}{
 author={Hutchinson, John E.},
 title={Fractals and self-similarity},
 journal={Indiana Univ. Math. J.},
 volume={30},
 date={1981},
 number={5},
 pages={713--747},
 doi={10.1512/iumj.1981.30.30055},
}

\bib{KKM22}{article}{
 author={Kinnunen, Juha},
 author={Kurki, Emma-Karoliina},
 author={Mudarra, Carlos},
 title={Characterizations of weak reverse H\"older inequalities on metric measure spaces},
 journal={Math. Z.},
 volume={301},
 date={2022},
 pages={2269--2290},
 doi={10.1007/s00209-022-02976-y},
}

\bib{Kit24}{article}{
 author={Kittle, Samuel},
 title={Absolutely continuous self-similar measures with exponential separation},
 journal={Ann. Sci. \'Ec. Norm. Sup\'er. (4)},
 volume={57},
 number={4},
 date={2024},
 pages={1191--1231},
 doi={10.24033/asens.2589},
 eprint={arXiv:2103.12684},
}

\bib{Kit25}{article}{
 author={Kittle, Samuel},
 title={Absolutely continuous Furstenberg measures},
 journal={Proc. Lond. Math. Soc. (3)},
 volume={131},
 number={1},
 date={2025},
 pages={e70072},
 doi={10.1112/plms.70072},
 eprint={arXiv:2305.05757},
}

\bib{KK25}{article}{
 author={Kittle, Samuel},
 author={Kogler, Constantin},
 title={On absolute continuity of inhomogeneous and contracting on average self-similar measures},
 date={2025},
 eprint={arXiv:2409.18936v4},
}

\bib{KK26}{article}{
 author={Kittle, Samuel},
 author={Kogler, Constantin},
 title={Entropy theory for random walks on Lie groups},
 journal={J. Theoret. Probab.},
 volume={39},
 date={2026},
 pages={Article~13},
 doi={10.1007/s10959-025-01468-z},
 eprint={arXiv:2501.05395},
}

\bib{KKsupp}{misc}{
 author={Kittle, Samuel},
 author={Kogler, Constantin},
 title={Supplementary material for On $L^p$ density of self-similar measures and Bernoulli convolutions},
 note={GitHub repository, \href{https://github.com/ckkogler/l1-to-lp-supplementary-material}{\nolinkurl{github.com/ckkogler/l1-to-lp-supplementary-material}}},
}

\bib{KogLean}{misc}{
 author={Kittle, Samuel},
 author={Kogler, Constantin},
 title={$L^p$ densities of self-similar measures and Bernoulli convolutions},
 note={Lean formalization, GitHub repository, \href{https://github.com/ckkogler/lpselfsimilar-lean}{\nolinkurl{github.com/ckkogler/lpselfsimilar-lean}}},
}

\bib{Mey63}{article}{
 author={Meyers, Norman G.},
 title={An $L^p$-estimate for the gradient of solutions of second order elliptic divergence equations},
 journal={Ann. Scuola Norm. Sup. Pisa (3)},
 volume={17},
 date={1963},
 number={3},
 pages={189--206},
}

\bib{PSS00}{article}{
 author={Peres, Yuval},
 author={Schlag, Wilhelm},
 author={Solomyak, Boris},
 title={Sixty years of Bernoulli convolutions},
 book={
  title={Fractal geometry and stochastics II},
  series={Progress in Probability},
  volume={46},
  publisher={Birkh\"auser},
  address={Basel},
 },
 date={2000},
 pages={39--65},
 doi={10.1007/978-3-0348-8380-1_2},
}

\bib{Shm19}{article}{
 author={Shmerkin, Pablo},
 title={On Furstenberg's intersection conjecture, self-similar measures, and the $L^q$ norms of convolutions},
 journal={Ann. of Math. (2)},
 volume={189},
 date={2019},
 number={2},
 pages={319--391},
 doi={10.4007/annals.2019.189.2.1},
}

\bib{SS16}{article}{
 author={Shmerkin, Pablo},
 author={Solomyak, Boris},
 title={Absolute continuity of self-similar measures, their projections and convolutions},
 journal={Trans. Amer. Math. Soc.},
 volume={368},
 date={2016},
 number={7},
 pages={5125--5151},
 doi={10.1090/tran6696},
}

\bib{Sol95}{article}{
 author={Solomyak, Boris},
 title={On the random series $\sum\pm\lambda^n$ (an Erd\H{o}s problem)},
 journal={Ann. of Math. (2)},
 volume={142},
 date={1995},
 number={3},
 pages={611--625},
 doi={10.2307/2118556},
}

\bib{Spa12}{article}{
 author={Spadaro, Emanuele},
 title={Nondoubling $A_\infty$ weights},
 journal={Adv. Calc. Var.},
 volume={5},
 date={2012},
 number={3},
 pages={345--354},
 doi={10.1515/acv.2011.019},
}

\bib{Tao11}{book}{
 author={Tao, Terence},
 title={An introduction to measure theory},
 series={Graduate Studies in Mathematics},
 volume={126},
 publisher={American Mathematical Society},
 address={Providence, RI},
 date={2011},
 doi={10.1090/gsm/126},
 note={\PrintDOI{10.1090/gsm/126}},
}

\bib{Var18}{article}{
 author={Varj\'u, P\'eter P.},
 title={Recent progress on Bernoulli convolutions},
 book={
  title={European Congress of Mathematics},
  publisher={European Mathematical Society},
  address={Z\"urich},
 },
 date={2018},
 pages={847--867},
 doi={10.4171/176-1/38},
 eprint={arXiv:1608.04210},
}

\bib{Var19}{article}{
 author={Varj\'u, P\'eter P.},
 title={Absolute continuity of Bernoulli convolutions for algebraic parameters},
 journal={J. Amer. Math. Soc.},
 volume={32},
 date={2019},
 number={2},
 pages={351--397},
 doi={10.1090/jams/916},
}

\bib{Yu22}{misc}{
 author={Yu, Han},
 title={Bernoulli convolutions with Garsia parameters in $(1,\sqrt{2}]$ have continuous density functions},
 date={2022},
 note={\href{https://arxiv.org/abs/2108.01008v2}{arXiv:2108.01008v2}},
}

\end{biblist}
\end{bibdiv}
\end{document}